\documentclass[11pt,leqno]{amsart}
\usepackage{amssymb,verbatim,enumerate,ifthen,times}
\usepackage[mathscr]{eucal}
\usepackage[utf8]{inputenc}
\usepackage[T1]{fontenc}
\usepackage{tabu}
\usepackage{dsfont}
\usepackage{accents}
\def\ubar#1{\underaccent{\bar}{#1}}

\def\N{\mathbb{N}}

\def\Q{\mathbb{Q}}
\def\R{\mathbb{R}}

\def\C{\mathrm{\bf C}}
\def\A{\mathrm{\bf A}}
\def\Z{\mathscr{Z}}

\def\diam{\mbox{\rm diam}}
\def\Ref(#1|#2){(#1\hspace{.3mm}|\hspace{.3mm}#2)}
\newtheorem{theorem}{Theorem}
\newtheorem*{theorem*}{Theorem}
\long\def\Thm#1#2{\ifthenelse{\equal{#1}{*}}{\begin{theorem*}#2\end{theorem*}}
             {\begin{theorem}\label{T#1}#2\end{theorem}}}
\newtheorem{Atheorem}{Theorem}

\def\thm#1{Theorem~\ref{T#1}}

\newtheorem{proposition}[theorem]{Proposition}
\newtheorem*{proposition*}{Proposition}
\long\def\Prp#1#2{\ifthenelse{\equal{#1}{*}}{\begin{proposition*}#2\end{proposition*}}
             {\begin{proposition}\label{P#1}#2\end{proposition}}}
\def\prp#1{Proposition~\ref{P#1}}

\newtheorem{corollary}[theorem]{Corollary}
\newtheorem*{corollary*}{Corollary}
\long\def\Cor#1#2{\ifthenelse{\equal{#1}{*}}{\begin{corollary*}#2\end{corollary*}}
             {\begin{corollary}\label{C#1}#2\end{corollary}}}
\def\cor#1{Corollary~\ref{C#1}}

\newtheorem{lemma}[theorem]{Lemma}
\newtheorem*{lemma*}{Lemma}
\long\def\Lem#1#2{\ifthenelse{\equal{#1}{*}}{\begin{lemma*}#2\end{lemma*}}
             {\begin{lemma}\label{L#1}#2\end{lemma}}}
\def\lem#1{Lemma~\ref{L#1}}

\theoremstyle{definition}
\newtheorem{definition}[theorem]{Definition}
\newtheorem*{definition*}{Definition}
\long\def\Defn#1#2{\ifthenelse{\equal{#1}{*}}{\begin{definition*}\rm #2\end{definition*}}
             {\begin{definition}\label{D#1}\rm #2\end{definition}}}

\newtheorem{remark}[theorem]{Remark}
\newtheorem*{remark*}{Remark}
\long\def\Rem#1#2{\ifthenelse{\equal{#1}{*}}{\begin{remark*}\rm #2\end{remark*}}
             {\begin{remark}\label{R#1}\rm #2\end{remark}}}

\newtheorem{example}{Example}
\newtheorem*{example*}{Example}
\long\def\Exa#1#2{\ifthenelse{\equal{#1}{*}}{\begin{example*}\rm #2\end{example*}}
             {\begin{example}\label{Ex#1}\rm #2\end{example}}}

\def\eq#1{{\rm(\ref{E#1})}}
\def\Eq#1#2{\ifthenelse{\equal{#1}{*}}
  {\begin{equation*}\begin{aligned}[]#2\end{aligned}\end{equation*}}
  {\begin{equation}\begin{aligned}\label{E#1}#2\end{aligned}\end{equation}}}

\begin{document}

\date{\today}

\title{On a functional equation related to two-variable Cauchy means}

\author[T. Kiss]{Tibor Kiss}
\author[Zs. P\'ales]{Zsolt P\'ales}
\address{
\textit{Address of the first author:} 
Institute of Mathematics, MTA-DE Research Group “Equations, Functions and Curves”,
Hungarian Academy of Sciences and University of Debrecen, P. O. Box 12, 4010 Debrecen, Hungary\newline
\textit{Address of the second author:} Institute of Mathematics, University of Debrecen, H-4032 Debrecen, Egyetem t\'er 1, Hungary
}
\email{\{kiss.tibor,pales\}@science.unideb.hu}


\thanks{The research of the first author was supported in part by the Hungarian Academy of Sciences and in part by the \'UNKP-17-Doctorand New National Excellence Program of the Ministry of Human Capacities. The research of the second author was supported by the 
Hungarian Scientific Research Fund (OTKA) Grant K111651.}

\begin{abstract}
In this paper, we are dealing with the solution of the functional equation
\Eq{*}{
\varphi\Big(\frac{x+y}2\Big)(f(x)-f(y))=F(x)-F(y),
}
concerning the unknown functions $\varphi,f$ and $F$ defined on a same open subinterval of the reals. Improving the previous results related to this topic, we describe the solution triplets $(\varphi,f,F)$ assuming only the continuity of $\varphi$. 

As an application, under natural conditions, we also solve the equality problem of two-variable Cauchy means and two-variable quasi-arithmetic means.
\end{abstract}

\maketitle

\section{Introduction}
Let $J$ be a nonempty open real interval. We say that a two-place function $M:J\times J\to\R$ is a \emph{two-variable mean on $J$} if it possesses the \emph{Mean Value Property}, that is, if
\Eq{*}{
\min(x,y)\leq M(x,y)\leq\max(x,y),\qquad(x,y\in J).}
If both of the above inequalities are sharp whenever $x\neq y$, then $M$ is called a \emph{strict mean}. To formulate the problem what we would like to deal with, we need the following two special classes of means.

\emph{The Class of Cauchy Means.} The Cauchy Mean Value Theorem states that if $G,H:J\to\R$ are differentiable functions then, for all distinct elements $x,y\in J$, there exists a point $u$ in the open interval determined by the points $x$ and $y$ such that the equality
\Eq{*}{
G'(u)\big(H(x)-H(y)\big)=H'(u)\big(G(x)-G(y)\big)
}
holds. If $H'$ does not vanish on $J$, then, by the Rolle Mean Value Theorem, the function $H$ must be injective on $J$, hence, in this case, the above equality is equivalent to
\Eq{*}{
\Big(\frac{G'}{H'}\Big)(u)=\frac{G(x)-G(y)}{H(x)-H(y)}.
}
Now, one can see easily that $u$ has to be unique if, in addition, the ratio $G'/H'$ is invertible on $J$. We note that, under the mentioned conditions, the invertibility of $G'/H'$ is equivalent to that it is continuous and strictly monotone on $J$. The precise formulation and the proof of this statement can be found as Remark 1 in the paper \cite{Mat04c}. To prove the essential part, the author uses the Darboux property of the derivative functions.

Motivated by the above observation, one can introduce the notion of Cauchy means as follows. We say that a mean $M:J\times J\to\R$ is a \emph{two-variable Cauchy mean} if there exist differentiable functions $G,H:J\to\R$ with $0\notin H'(J)$ such that the function $G'/H'$ is invertible on $J$ and that, for all $x,y\in J$, we have
\Eq{CM}{
M(x,y)=\Big(\dfrac{G'}{H'}\Big)^{-1}\Big(\dfrac{G(x)-G(y)}{H(x)-H(y)}\Big)\quad\mbox{if }x\neq y\qquad\mbox{and}\qquad M(x,y)=x\quad\mbox{if }x=y.}
In this case, we denote $M$ by $\C_{G, H}$, where $G$ and $H$ will be called the \emph{generator functions of the Cauchy mean}. There are many papers dealing with this class of means. In particular, the equality and homogeneity problem of Cauchy means was completely solved by Losonczi in \cite{Los00a,Los03a} and \cite{Los02a}, respectively. For the solution of the equality problem in \cite{Los00a}, 7th-order differentiability was assumed. This regularity condition was reduced to first-order differentiability by Matkowski in \cite{Mat04c}. The comparison problem of Cauchy means was also studied by Losonczi \cite{Los02b}. The so-called invariance equation for Cauchy means was solved by Berrone in \cite{Ber05}. Characterization type results were obtained by Berrone \cite{Ber16} and by Berrone--Moro \cite{BerMor00}.

\emph{The Class of Quasi-Arithmetic Means.} We say that a mean $M:J\times J\to\R$ is a \emph{two-variable quasi-arithmetic mean} if one can find a continuous strictly monotone function $\Phi:J\to\R$ such that, for all $x,y\in J$, we have
\Eq{*}{
M(x,y)=\Phi^{-1}\Big(\frac{\Phi(x)+\Phi(y)}2\Big)=:\A_\Phi(x,y).
}
The continuity and the strict monotonicity of $\Phi$ provides that $\Phi$ is invertible and that $\Phi(J)$ is an open subinterval of $\R$. Hence the above expression is well-formulated and it indeed defines a strict mean on $J$. Similarly to the previous part, $\Phi$ will be called the \emph{generator function of the quasi-arithmetic mean}.

In this paper we are going to focus on the \emph{equality problem of two-variable Cauchy means and two-variable quasi-arithmetic means}. It is easy to check that the equality problem $\C_{G,H}=\A_\Phi$ on $J$ can be rewritten as
\Eq{-1.5}{
\varphi\Big(\frac{u+v}2\Big)\big(f(u)-f(v)\big)=F(u)-F(v),\qquad(u,v\in \Phi(J)),
}
where we define
\Eq{-1.4}{
\varphi:=\frac{G'}{H'}\circ\Phi^{-1},\qquad
f:=H\circ\Phi^{-1},\qquad\mbox{and}\qquad
F:=G\circ\Phi^{-1}.
}

Motivated by this reformulation, to solve the equality problem, it is enough to describe the solution triplets $(\varphi,f,F)$ of the functional equation
\Eq{-1}{
\varphi\Big(\frac{x+y}2\Big)\big(f(x)-f(y)\big)=F(x)-F(y),\qquad(x,y\in I),
}
where $I$ stands for a nonempty open subinterval of $\R$.

The functional equations having the form \eq{-1} have a rich literature and it was investigated by several authors. The firs remarkable result, from the year 1985, is due to J. Aczél \cite{Acz85a}, who solved \eq{-1} assuming that $f$ is the identity function. This particular case was also dealt with by Volkmann \cite{Vol84c}. Without any regularity assumptions for the unknown functions, he proved that the pair $(\varphi,F)$ solves the equation if and only if $F$ is a polynomial of degree at most two and $\varphi$ is its derivative. Independently, in 1979, Sh. Haruki investigated a pexiderized version of the equation of Aczél, where $f$ was still the identity, and he proved the same result, again without any regularity conditions. The details of Haruki's work can be found in \cite{Har79}. The solutions of the equation in \eq{-1} was first described in 2016, by Z. M. Balogh, O. O. Ibrogimov and B. S. Mityagin in the paper \cite{BalIbrMit16}. In their approach, the functions $f$ and $g$ was supposed to be three-times differentiable. Two years later, R. Łukasik improved this result in \cite{Luk18}, by showing that the first-order differentiability of $f$ and $g$ is sufficient to obtain the same solutions. Now we present a different approach, which allows us to treat and solve \eq{-1} solely under the continuity of $\varphi$. 

The functional equation \eq{-1} seems to be a particular case of the functional equation 
\Eq{*}{
  \varphi(x+y)=\frac{F(x)G(y)+H(x)L(y)}{m(x)+n(y)}
}
which was studied by Lundberg in a series of papers \cite{Lun99}, \cite{Lun01}. However, the solutions therein have indirect forms and we could not find way to deduce our results from those of \cite{Lun99} and \cite{Lun01}. Another problem, which was formulated by P. K. Sahoo and T. Riedel in their book \cite[Section 2.7]{SahRie98}, remains still open: \emph{without any regularity conditions concerning the members of $(\varphi,f,\Psi,F)$}, solve the functional equation
\Eq{*}{
\varphi\Big(\frac{x+y}2\Big)\big(f(x)-f(y)\big)=
\Psi\Big(\frac{x+y}2\Big)\big(F(x)-F(y)\big),\qquad(x,y\in\R).
}

One of our main observation (which will turn out in the next section) is that equation \eq{-1} has a strong connection to the functional equation
\Eq{-.5}{
	\varphi\Big(\frac{x+y}2\Big)(f(x)+f(y))=F(x)+F(y),\qquad(x,y\in I),
}
which has been studied and solved by the authors in the paper under minimal regularity assumption in \cite{KisPal18}.

\section{Connection among the functional equations \eq{-1} and \eq{-.5}}

We say that a subinterval $J\subseteq I$ is \emph{proper} if it has at least two elements and is different from $I$. A function $g:I\to\R$ will be called \emph{locally non-constant on $I$} if there is no proper subinterval of $I$, where $g$ is constant. In other words, $g$ is locally non-constant on $I$ if and only if, for all $\lambda\in\R$, the set $g^{-1}(\{\lambda\})$ has an empty interior in the domain $I$.

First we refer to the result about the solutions of equation \eq{-.5}. Observe, that the substitution $x=y$ in the equation \eq{-.5} immediately yields that $F=\varphi \cdot f$ holds on $I$. Therefore, in view of \cite[Theorem 11]{KisPal18}, we have the following statement.

\Thm{m1}{
Let $(\varphi,f,F):I\to\R^3$ such that $\varphi$ and $f$ are continuous on $I$. 
Then the triplet $(\varphi,f,F)$ solves functional equation \eq{-.5} if and only if either
\begin{enumerate}
\item[(i)] there exists an interval $J\subseteq I$ such that $f(x)=0$ for all $x\in I\setminus J$, the function $\varphi$ is constant on $\frac12(J+I)$, and $F=\varphi\cdot f$,
\end{enumerate}
or
\begin{enumerate}
\item[(ii)] $f$ is nowhere zero on $I$ and there exist constants $a,b,c,d,\gamma\in\R$ with $ad\neq bc$ such that
\Eq{m1a}{
\begin{array}{c}
 F(x)=c\sin(\sqrt{-\gamma}x)+d\cos(\sqrt{-\gamma}x), \\[1mm]
 f(x)=a\sin(\sqrt{-\gamma}x)+b\cos(\sqrt{-\gamma}x),
\end{array}
\qquad 
\varphi(x)=\dfrac{c\sin(\sqrt{-\gamma}x)+d\cos(\sqrt{-\gamma}x)}{a\sin(\sqrt{-\gamma}x)+b\cos(\sqrt{-\gamma} x)} 
\quad \mbox{ if } \, \gamma<0,
}
\Eq{m1b}{
\begin{array}{c}
 F(x)=cx+d, \\[1mm]
 f(x)=ax+b,
\end{array}
\qquad \varphi(x)=\dfrac{cx+d}{ax+b}\quad \mbox{ if } \, \gamma=0,
}
\Eq{m1c}{
\begin{array}{c}
 F(x)=c\sinh(\sqrt{\gamma}x)+d\cosh(\sqrt{\gamma}x), \\[1mm]
 f(x)=a\sinh(\sqrt{\gamma}x)+b\cosh(\sqrt{\gamma}x),
\end{array}
\qquad  
\varphi(x)=\dfrac{c\sinh(\sqrt{\gamma}x)+d\cosh(\sqrt{\gamma}x)}{a\sinh(\sqrt{\gamma}x)+b\cosh(\sqrt{\gamma}x)}
\quad \mbox{ if } \, \gamma>0
}
\end{enumerate}
holds for all $x\in I$.
}

It is worth noticing that the condition \emph{(ii)} of the above theorem could be formulated in the following equivalent way:
\begin{enumerate}
 \item[(ii)'] $f$ is nowhere zero on $I$, there exists a constant $\gamma\in\R$ such that $f$ and $F$ are linearly independent solutions of the second-order linear homogeneous differential equation $Y''=\gamma Y$, and $\varphi=F/f$ on $I$.
\end{enumerate}

In fact, in the paper \cite[Theorem 11]{KisPal18}, instead of the continuity of $f$, a weaker regularity assumption was made. However, the following consequence of \thm{m1}, in which the regularity assumptions for $f$ are completely eliminated, will be sufficient for our purposes.

\Cor{m1}{Let $(\varphi,f,F):I\to\R^3$ such that $\varphi$ is continuous and locally non-constant on $I$. 
Then the triplet $(\varphi,f,F)$ solves functional equation \eq{-.5} if and only if $F=\varphi \cdot f$ and either
$f=0$ and $\varphi$ is arbitrary on $I$, or the alternative (ii) of \thm{m1} holds. Consequently, $f$ is an infinitely many times differentiable function.
}

\begin{proof}
First we show that the assumptions made on $\varphi$ and \eq{-.5} imply that $f$ is continuous on $I$. To do this, let $x_0\in I$ 
be arbitrarily fixed. Then there exists $r>0$ such that $[x_0-2r,x_0+2r]\subseteq I$. The function $\varphi$ is non-constant on the interval $[x_0-r,x_0+r]$, consequently, there exists $u_0\in[x_0-r,x_0+r]$ such that $\varphi(x_0)\neq\varphi(u_0)$. Let $y_0:=2u_0-x_0$. Then $y_0\in [x_0-2r,x_0+2r]$, the function 
$x\mapsto\varphi(\frac{x+y_0}{2})-\varphi(x)$ is continuous on $I$, and, by the choice of $u_0$, it is different from zero at the point $x:=x_0$. Therefore, there exists $\delta>0$ such that this function is different from zero also on the entire interval $]x_0-\delta,x_0+\delta[\,\subseteq I$. 
Using this, and the equality $F=\varphi \cdot f$, equation \eq{-.5} directly implies that
\Eq{*}{
f(x)
=f(y_0)\frac{\varphi(y_0)-\varphi(\frac{x+y_0}{2})}{\varphi(\frac{x+y_0}{2})-\varphi(x)}
,\qquad(x\in\,]x_0-\delta,x_0+\delta[\,).
}
Thus $f$ coincides with a continuous function on the neighborhood $]x_0-\delta,x_0+\delta[\,$ of $x_0$, 
which means that $f$ has to be continuous at $x_0$. Because $x_0$ was arbitrarily chosen, it follows that $f$ is continuous on $I$.

Thus, by \thm{m1}, we have the alternatives \emph{(i)} and \emph{(ii)} for the solutions of \eq{-.5}. The function $f$ is obviously infinitely many times differentiable, provided that we have the case \emph{(ii)}. On the other hand, having \emph{(i)} and using that $\varphi$ is locally non-constant, it follows that the subinterval $J$ in alternative \emph{(i)} must be empty. Hence, in this case, $f$ is identically zero on $I$, which finishes the proof. 
\end{proof}

In order to describe the connection between functional equations \eq{-1} and \eq{-.5}, let us introduce the following notations. For a positive number $h$ and a function $f:I\to\R$, define the subinterval interval $I_h$ of $I$ and the function $\delta_h f:I_h\to\R$ by
\Eq{*}{
  I_h:=(I-h)\cap(I+h) \qquad\mbox{and}\qquad \delta_hf(x):=f(x+h)-f(x-h), \qquad(x\in I_h).
}
With these notations, we have the following basic result, which derives a functional equation from \eq{-1} for $\varphi$ and $\delta_hf$, furthermore eliminates the function $F$.

\Thm{0}{
Assume that the triplet of functions $(\varphi,f,F):I\to\R^3$ solves functional equation \eq{-1}. Then, for all $h>0$, we have
\Eq{h}{
\varphi\Big(\frac{x+y}2\Big)(\delta_hf(x)+\delta_hf(y))=\varphi(x)\delta_hf(x)+\varphi(y)\delta_hf(y),\qquad(x,y\in I_h),
}
that is, the triplet $(\varphi,\delta_hf,\varphi\cdot\delta_hf)$ solves functional equation \eq{-.5} on $I_h$. Furthermore, if $\varphi$ is continuous and locally non-constant on $I$, then $f$ is infinitely many times differentiable on $I$ and
\Eq{f'}{
\varphi\Big(\frac{x+y}2\Big)(f'(x)+f'(y))=\varphi(x)f'(x)+\varphi(y)f'(y),
\qquad(x,y\in I),
}
that is, the triplet $(\varphi,f',\varphi\cdot f')$ satisfies functional equation \eq{-.5} on $I$.}

\begin{proof}
Assume that equation \eq{-1} holds.
Let $h>0$ and $x,y\in I_h$. Then $x\pm h$ and $y\pm h$ belong to $I$. Substituting the pairs $(x-h,x+h)$, $(x+h,y-h)$, $(y-h,y+h)$ and $(y+h,x-h)$ into \eq{-1}, we get the following four equalities:
\Eq{*}{
  \varphi(x)(f(x-h)-f(x+h))&=F(x-h)-F(x+h),\\
  \varphi\Big(\frac{x+y}2\Big)(f(x+h)-f(y-h))&=F(x+h)-F(y-h),\\
  \varphi(y)(f(y-h)-f(y+h))&=F(y-h)-F(y+h),\\
  \varphi\Big(\frac{x+y}2\Big)(f(y+h)-f(x-h))&=F(y+h)-F(x-h).
}
Adding up the above equations side by side, the equality in \eq{h} follows immediately.

Assume now that $\varphi$ is continuous and locally non-constant on $I$.
Applying \cor{m1} for the function $\delta_hf$ (instead of $f$), it follows, for all positive real number $h$, that the function $\delta_h f$ is infinitely many times differentiable on the subinterval $I_h\subseteq I$. Thus, based on the celebrated result of N. G. de Bruijn \cite{Bru51,Bru52}, $f$ can be written as $f_0+g$, where $f_0:I\to\R$ is infinitely many times differentiable and $g:\R\to\R$ is additive. Without loss of generality, we may assume that $g$ vanishes on the set of the rationals or, equivalently, that $f=f_0$ on the set $I\cap\Q$. This immediately implies that $(\varphi,f_0,F)$ solves the equation \eq{-1} on the intersection $I\cap\Q$.

Now, indirectly, assume that $g$ is not identically zero. Then there exists $x_0\in I\setminus\Q$ such that $g(x_0)\neq 0$. Let $(x_n)\subseteq I\cap\Q$ be any sequence such that $x_n\to x_0$ as $n\to\infty$. Applying the equation \eq{-1} for the triplet $(\varphi,f_0,F)$, for any member of $(x_n)$, and for any point $y\in I\cap\Q$, we get that
\Eq{*}{
F(x_n)=\varphi\Big(\frac{x_n+y}2\Big)(f_0(x_n)-f_0(y))+F(y),\qquad(n\in\N).
}
The continuity of $\varphi$ and $f_0$ provides that the limit of the left hand side exists as $n\to\infty$. Denoting $\lim\nolimits_{n\to\infty}F(x_n)$ by $\lambda$, using the decomposition of $f$, and the equality $f=f_0$ on $I\cap\Q$, we get that
\Eq{*}{
\lambda
&=\varphi\Big(\frac{x_0+y}2\Big)(f_0(x_0)-f_0(y))+F(y)
=\varphi\Big(\frac{x_0+y}2\Big)(f(x_0)-f(y))+F(y)-\varphi\Big(\frac{x_0+y}2\Big)g(x_0).
}
On the very right hand side, we can apply that $(\varphi,f,F)$ solves equation \eq{-1} on $I$, hence
\Eq{*}{
\lambda=F(x_0)-\varphi\Big(\frac{x_0+y}2\Big)g(x_0),\qquad(y\in I\cap\Q).
}
By our assumption, $g(x_0)\neq 0$, thus this equation is equivalent to
\Eq{*}{
\varphi\Big(\frac{x_0+y}2\Big)=\frac{F(x_0)-\lambda}{g(x_0)},\qquad(y\in I\cap\Q),
}
which means that the function $\varphi$ is constant on the dense subset $\frac12(x_0+(I\cap\Q))$ of the interval $\frac12(x_0+I)$. Thus, by its continuity, $\varphi$ must be constant on the subinterval $\frac12(x_0+I)\subseteq I$, which contradicts the fact that it is locally non-constant on $I$. Consequently, $g$ must vanish on $I$, and hence also on $\R$. Therefore $f=f_0$ on $I$, where $f_0$ was an infinitely many times differentiable function.

Finally, we prove that \eq{f'} is also valid. Let $x,y\in I$ be arbitrary. Then, for small positive $h$, we have that $x,y\in I_h$ and hence \eq{h} holds for $x,y$ and for small positive $h$. Dividing equation \eq{h} by $2h$ and taking the limit $h\to 0^+$, we get that \eq{f'} is satisfied.
\end{proof}

In order to be able to apply \thm{0} to establish the solutions of \eq{-1}, we have to distinguish two main cases concerning $\varphi$. First we will deal with the case when $\varphi$ is constant on a proper subinterval of $I$ and then with the case when $\varphi$ is \emph{locally non-constant}. The second part of the result of \thm{0} applies in the latter case. The investigation of the first case requires a detailed analysis, thus, firstly, we turn to this part.

\section{Preliminary results}

In order to treat \eq{-1}, first we investigate a special case, more precisely, we solve the much simple functional equation
\Eq{0}{
\varphi\Big(\frac{x+y}{2}\Big)(f(x)-f(y))=0,\qquad(x,y\in I),
}
which contains the two unknown functions $\varphi:I\to\R$ and $f:I\to\R$. We note that a functional equation having the form as in \eq{0} can be derived from \eq{-1} assuming that $F$ is an affine transformation of the function $f$, that is, there exists $A,B\in\R$ such that $F=Af+B$ on $I$.

For the brevity, we will frequently use the notations $\alpha:=\inf I$ and $\beta:=\sup I$. If $S\subseteq I$, then the complementary set $I\setminus S$ will be simply denoted by $S^c$. For a given function $g:I\to\R$, denote the set $g^{-1}(\{0\})$ by $\Z_g$. For two subsets $S,P\subseteq\R$, let us denote the set $(2P-S)\cap I$ by $\Ref(S|P)$. That is, $\Ref(S|P)$ consists of those elements of $I$ that are reflections of an element of $S$ with respect to some element of $P$. 

\Prp{HK}{
Let $P$ and $S$ be arbitrary subsets of $\R$.
\begin{enumerate}[(a)]\itemsep=1mm
\item\label{a} If $P_1\subseteq P_2\subseteq\R$ and $S_1\subseteq S_2\subseteq\R$, then $\Ref(S|P_1)\subseteq\Ref(S|P_2)$ 
and $\Ref(S_1|P)\subseteq\Ref(S_2|P)$.
\item If at least one of the sets $P$ or $S$ is open, then $\Ref(S|P)$ is also open.
\item\label{b} If $P,S\subseteq\R$ are intervals, then $\Ref(S|P)$ is also an interval, furthermore we have
\Eq{ip}{
\inf\Ref(S|P)=\max(\alpha,2\inf P-\sup S)
\qquad\mbox{and}\qquad
\sup\Ref(S|P)=\min(2\sup P-\inf S,\beta).
}
\item\label{c} Consequently, for a given point $p\in I$, the set $\Ref(I|p)$ is the maximal open subinterval contained in $I$, which is symmetric with respect to $p$. Furthermore, by the equality \eq{ip}, it follows that $\Ref(I|p)$ is bounded unless $I=\R$.
\end{enumerate}
}

\begin{proof}
The statements \emph{(a)}, \emph{(b)}, \emph{(c)}, and \emph{(d)} are direct consequences of the definition.
\end{proof}

\Lem{0}{
Let $(\varphi,f)$ be a solution of equation \eq{0} and $p\in I$. Then $f(x)=f(p)$ holds for all 
$x\in \Ref(p|{\Z_\varphi}^c)$.}

\begin{proof}
If $\Ref(p|{\Z_\varphi}^c)$ is empty, then the statement is obvious. Therefore we may assume that $\Ref(p|{\Z_\varphi}^c)$ is nonempty 
and let $x\in \Ref(p|{\Z_\varphi}^c)$ be arbitrary. Then there exists $u\in {\Z_\varphi}^c$ such that $\frac{x+p}2=u$ holds. Applying 
equation \eq{0} for $x$ and $p$ and using that $\varphi(u)\neq 0$, we get $f(x)=f(p)$.
\end{proof}

\Prp{1}{
Let $(\varphi,f)$ be a solution of \eq{0} such that $\Z_\varphi$ is closed in $I$ and has an empty interior. Then $f$ is constant on $I$.
}

\begin{proof}
The assumptions concerning $\Z_\varphi$ provide that ${\Z_\varphi}^c$ is a nonempty open subset of $I$. Therefore it can be written as a 
union 
of its components. Let $K\subseteq{\Z_\varphi}^c$ be any component and $x,y\in K$ be arbitrary. Then $\frac{x+y}2\in K$, that is 
$\varphi(\frac{x+y}{2})\neq0$, hence, based on equation \eq{0}, we get that $f(x)=f(y)$.

In the next step we show that $f$ takes the same value on any two components of ${\Z_\varphi}^c$. There is nothing to prove if 
${\Z_\varphi}^c$ has only one component. Thus we may assume that this is not the case and let $K$ and $L$ be different components of 
${\Z_\varphi}^c$. By our assumption, the interior of $\Z_\varphi$ is empty, hence the intersection $\frac{1}{2}(K+L)\cap({\Z_\varphi}^c)$ 
cannot be empty. Let $u$ be any point of the intersection. Then $\varphi(u)\neq 0$, furthermore, there exist $x\in K$ and $y\in L$ such 
that $\frac{x+y}{2}=u$. Using equation \eq{0} for the points $x$ and $y$, we obtain that $f(x)=f(y)$. Consequently, there exists 
$\lambda\in\R$ such that $f(x)=\lambda$ whenever $x\in {\Z_\varphi}^c$.

Finally we show that the restriction $f|_{\Z_\varphi}$ is also identically $\lambda$. If $\Z_\varphi$ is empty, then we are done. If 
not then let $p\in \Z_\varphi$ be arbitrarily fixed. The assumptions concerning $\Z_\varphi$ yield that there is an open interval 
$U\subseteq {\Z_\varphi}^c$ such that $J:=\Ref(p|U)\subseteq \Ref(p|{\Z_\varphi}^c)$ is contained in $I$. Then, by \lem{0} we 
have $f|_J\equiv f(p)$. On the other hand, $J$ is a nonempty interval in $I$, which implies that the intersection $J\cap{\Z_\varphi}^c$ 
cannot be empty. Consequently, $f(p)=\lambda$ must hold. 
\end{proof}

In the following proposition, we establish the solutions of equation \eq{0} supposing that the interior of $\Z_\varphi$ is nonempty. In this case, the interior of $\Z_\varphi$ is the union if its components which are maximal open subintervals of the interior of $\Z_\varphi$.
Provided that $\Z_\varphi$ is closed, the closure of such intervals are maximal closed subintervals of $\Z_\varphi$.

\Prp{2}{
Let $(\varphi,f)$ be a solution of \eq{0} such that $\Z_\varphi$ is a closed subset of $I$ with a nonempty interior. Then, for any 
maximal closed subinterval $J\subseteq\Z_\varphi$, there exist constants $\lambda,\mu\in\R$ such that
\Eq{*}{
\Ref(I|\inf J)\subseteq\Z_{f-\lambda}
\qquad\mbox{and}\qquad
\Ref(I|\sup J)\subseteq\Z_{f-\mu}.}}

\begin{proof} Let $J$ be any nonempty closed subinterval of $\Z_\varphi$. If $J=I$, then $\varphi$ is identically zero on $I$ and, 
by $\Ref(I|\inf J)=\Ref(I|\sup J)=\emptyset$, the inclusions in the statement above hold for any constants $\lambda$ and $\mu$. Thus, 
we may assume that at least one of the endpoints of $J$ belongs to $I$, say $p:=\sup J\in I$. Then $\Ref(I|p)$ is nonempty. 
Let $a:=\inf \Ref(I|p)$ and $b:=\sup \Ref(I|p)$.

In the first step we are going to show that $f$ is constant on $]a,p]$. Let $x<y$ in $]a,p]$ be arbitrary. Then $p<\frac{x+\beta}{2}$ 
and, due to the maximality of $J$ and the openness of ${\Z_\varphi}^c$, one can find a nonempty open and bounded subinterval $U$ contained 
in ${\Z_\varphi}^c\cap\,\big]p,\frac{x+\beta}{2}\big[$.

Let $r:=2|U|$ and consider first the case $y-x<r$. Then there is $u<v$ in $U$ such that $v-u=\frac12(y-x)$, that is $2v-y=2u-x=:w$. A short 
calculation yields that $w\in I$, which means that the intervals $\Ref(y|U)$ and $\Ref(x|U)$ have a common point in $I$. By In view of 
\lem{0}, $f(x)=f(y)$ follows.

If $y-x\geq r$, then there exist $n\in\N$ and $t_0<\dots<t_n$ in $[x,y]$ such that $t_0=x$, $t_n=y$ and $t_i-t_{i-1}<r$ for all 
$i\in\{1,\dots,n\}$. Based on the previous argument, $f(t_{i-1})=f(t_i)$ follows for all $i\in\{1,\dots,n\}$, particularly, we have again 
$f(x)=f(y)$. Thus, there exists $\mu\in\R$ such that $f(x)=\mu$ if $x\in\,]a,p]$.

Now we show that $f$ takes $\mu$ also on $]p,b[$. Let $x\in\,]p,b[$ be any point. Then, by the maximality of $J$, there exists 
$u\in{\Z_\varphi}^c$ with $p<u<x$. Then $2u-x$ belongs to $]a,p]$. Applying equation \eq{0} and using $\varphi(u)\neq 0$, 
we immediately get that $f(x)=f(2u-x)=\mu$.
\end{proof}

To formulate our next theorem, for a given set $S\subseteq I$, we introduce the notations
\Eq{*}{
S_*:=\{x\in I\mid x<\inf S\}
\qquad\mbox{and}\qquad
S^*:=\{x\in I\mid \sup S<x\}.
}
We note that, regardless of the topological properties of $S$, the sets $S_*$ and $S^*$ are always open in $I$ and that $S_*=S^*=I$ if and 
only if $S$ is empty.

\Lem{s}{Let $S\subseteq I$ be nonempty, $s_*:=\inf\frac12(S+I)$ and $s^*:=\sup\frac12(S+I)$. Then we have
\Eq{seq}{
S_*=\Ref(I|s_*)
\qquad\mbox{or}\qquad
S^*=\Ref(I|s^*),
}
provided that $-\infty<s_*$ or $s^*<+\infty$, respectively.
}

\begin{proof}
We prove only the first identity in \eq{seq}, the second one can be shown similarly. Assume that $-\infty<s_*$ holds. Then, $s_*$ is finite and by its definition,
\Eq{s*}{
  s_*=\frac{\inf S+\alpha}{2}<\frac{\beta+\alpha}{2}.
}
Then $\alpha$ is also finite and
\Eq{*}{
  S_*=\,]\alpha,\inf S[=\,]\alpha,2s_*-\alpha[.
}
On the other hand, by formula \eq{ip} and the inequality in \eq{s*},
\Eq{*}{
  (I|s_*)=\,]\max(\alpha,2s_*-\beta),\min(2s_*-\alpha,\beta)[=\,]\alpha,2s_*-\alpha[.
}
which finishes the proof.
\end{proof}

\Thm{1}{
Let $\varphi,f:I\to\R$ such that $\Z_\varphi$ is closed. Then the pair $(\varphi,f)$ is a solution of equation \eq{0} if and only if
either $f$ is constant on $I$ and $\varphi:I\to\R$ is any function or there exist a nonempty closed interval $K\subseteq I$ and 
constants $\lambda_*,\lambda^*\in\R$ such that the inclusions
\Eq{Sol0}{
K_*\subseteq\Z_{f-\lambda_*},
\qquad
K^*\subseteq\Z_{f-\lambda^*}
\qquad\mbox{and}\qquad
\tfrac12(I+K)\subseteq\Z_\varphi,
}
are satisfied.}

\begin{proof}
If $f$ is constant on $I$, then, obviously, for any function $\varphi:I\to\R$, the pair $(\varphi,f)$ is a solution of \eq{0}.
Let $K\subseteq I$ be a nonempty closed interval, $\lambda_*,\lambda^*\in\R$ such that \eq{Sol0} holds, finally let $x,y\in I$
be arbitrarily fixed. If either $x,y\in K_*$ or $x,y\in K^*$, then either $f(x)=f(y)=\lambda_*$ or $f(x)=f(y)=\lambda^*$, 
respectively, which imply the equality in \eq{0}. In the remaining case, we have that $\inf K\leq\max(x,y)$ and $\min(x,y)\leq\sup K$,
therefore
\Eq{*}{
\inf(K+I)=\inf K+\alpha<\max(x,y)+\min(x,y)<\beta+\sup K=\sup(K+I),
}
where, obviously, $\max(x,y)+\min(x,y)=x+y$.
These mean that $\frac{x+y}{2}\in\frac12(K+I)\subseteq\Z_\varphi$, that is, $\varphi\big(\frac{x+y}{2}\big)=0$, whence the equality in 
\eq{0} follows again.  

Conversely, assume that $(\varphi,f)$ is a solution of equation \eq{0} such that $f$ is non-constant on $I$. If $\Z_\varphi=I$ then the statement of the theorem holds with $K=I$. (Then $K_*$ and $K^*$ are empty and $\lambda_*$, $\lambda^*$ can be arbitrary.) Thus, in the rest of the proof, we may assume that $\Z_\varphi$ is a proper subset of $I$.

Let $J:=\,]\alpha,a[$ and $L:=\,]b,\beta[\,$, where
\Eq{*}{
a:=\sup\{t\geq\alpha\mid f\text{ is constant on }]\alpha,t[\}
\qquad\text{and}\qquad
b:=\inf\{t\leq\beta\mid f\text{ is constant on }]t,\beta[\}.
}
Then $J$ and $L$ are maximal intervals with respect to the property that $f$ is constant on them. Let $K:=(J\cup L)^c$. By \prp{1}, the interior of $\Z_\varphi$ cannot be empty. Thus, $\Z_\varphi$ contains a maximal closed subinterval $J_0$ which cannot be equal to $I$ because $\Z_\varphi$ is a proper subset of $I$. Thus, one of the endpoints of $J_0$, say $p$ is in $I$. In view of \prp{2}, it follows that $f$ is constant on $(I|p)$, which is a nonempty interval contained either in $J$ or in $L$. Therefore at least one of the intervals $J$ or $L$ is nonempty, consequently, $K$ is different from $I$. By our assumption concerning $f$, the intervals $J$ and $L$ are disjoint and different from $I$, which implies that $K$ is nonempty. We have obtained that $K$ is a nonempty closed proper subinterval of $I$. Furthermore, we also have $K_*=J$ and $K^*=L$, thus there exist $\lambda_*,\lambda^*\in\R$ such that $K_*\subseteq\Z_{f-\lambda_*}$ and $K^*\subseteq\Z_{f-\lambda^*}$. Thus, the first two inclusions in \eq{Sol0} hold and we only need to prove that $\varphi$ is identically zero on $G:=\frac12(K+I)$.

By $K\neq I$, the interval $G$ strictly contains $K$, thus $G$ and $J\cup L$ have common points. This yields that $f$ cannot be 
constant on $G$. Restricting the equation \eq{0} to $G$, by \prp{1}, we get that the interior of $G\cap\Z_\varphi$ is nonempty. Thus 
there exists a maximal nonempty open subinterval $U_0\subseteq G$ such that $\varphi$ is identically zero on $U_0$. 

To complete the proof, it suffices to show that $U_0=G$. Assume, to the contrary, that this is not the case. 
Then we may also assume that $s_*:=\inf G<p:=\inf U_0$ holds. Let $U\subseteq I$ be a maximal closed subinterval such that 
$U_0\subseteq U$ and $\varphi$ is identically zero on $U$. Then $\inf U$ must be also equal to $p$. Applying \prp{2} for the maximal 
closed subinterval $U\subseteq\Z_\varphi$, the function $f$ is constant on $\Ref(I|p)$. Then exactly one of the identities 
$\inf\Ref(I|p)=\alpha$ or $\sup\Ref(I|p)=\beta$ can hold. By the maximality property of $J$ and $L$, we have either 
$\Ref(I|p)\subseteq J$ or $\Ref(I|p)\subseteq L$, respectively. Because of the symmetry, we may assume that $\Ref(I|p)\subseteq J$. Then $J$ is nonempty, 
which implies that $\alpha<a$ must hold. To finish the proof, we distinguish two cases.

\emph{Case 1.} If $\alpha=s_*$ then, having the relation $\alpha<a$, it follows that $I$ is not bounded from below, that is, 
$\alpha=s_*=-\infty=\inf\Ref(I|p)$. Then there exists $x\in\Ref(I|p)$ such that $x<2p-a$, which yields that $a<2p-x\in\Ref(I|p)\subseteq J$. This 
contradicts that $\sup J=a$.

\emph{Case 2.} Assume now that $\alpha<s_*$. Then, using \lem{s} and that $\inf\Ref(I|s_*)=\inf\Ref(I|p)=\alpha$, $s_*<p$, we obtain that
\Eq{*}{
J=K_*=\Ref(I|s_*)\subsetneq\Ref(I|p)\subseteq K_*=J,}
which is a contradiction again.

Consequently, $\frac12(K+I)=G=U_0\subseteq\Z_\varphi$. 
\end{proof}

\Rem{rem1}{One can easily see that the closedness of the set $\Z_\varphi$ was only used to establish the necessity of the statement of \thm{1}. The sufficiency holds without any assumption on $\varphi$.}

\section{Solutions of \eq{-1} assuming that $\varphi$ is constant on a subinterval of $I$}

In this section we are going to describe the solutions of the functional equation \eq{-1} assuming that the function $\varphi$ is constant on a proper subinterval of $I$. 
The next proposition clarifies the main connection between the equations \eq{-1} and \eq{0}.

\Prp{3}{Let $(\varphi,f,F)$ be a solution of equation \eq{-1} on the interval $I$, $A\in\R$, and $J\subseteq I$ be a subinterval. Then the pair $(\varphi_A,f)$ solves the equation \eq{0} on $J$ if and only if there exists $B\in\R$ such that $F=Af+B$ on $I$.
}

\begin{proof}
Assume that $(\varphi_A,f)$ solves \eq{0} on the interval $J$, that is, we have
\Eq{*}{
\varphi_A\Big(\frac{x+y}{2}\Big)(f(x)-f(y))=
\Big(\varphi\Big(\frac{x+y}{2}\Big)-A\Big)(f(x)-f(y))=0,\qquad(x,y\in J).
}
On the other hand, particularly, equation \eq{-1} holds on $J$. Subtracting the above equation from \eq{-1} side by side, we get that
\Eq{*}{
F(x)-Af(x)=F(y)-Af(y),\qquad(x,y\in J).
}
Consequently, there exists $B\in\R$ such that $F(x)-Af(x)=B$ for all $x\in J$, thus $F=Af+B$ holds on $J$.
	
Conversely, assume that there are constants $A,B\in\R$ such that $F=Af+B$ on the interval $J$. Using this and that $(\varphi,f,F)$ is a solution of equation \eq{-1} on $J$, we obtain that
\Eq{*}{
\varphi\Big(\frac{x+y}{2}\Big)(f(x)-f(y))&=F(x)-F(y)=Af(x)+B-Af(y)-B=A(f(x)-f(y))
}
for all $x,y\in J$. This means that $(\varphi_A,f)$ solves \eq{0} on $J$, which finishes the proof.
\end{proof}

The following consequence of \prp{3} will play a key role in the proof of \thm{2} below.

\Cor{ext}{
Let $(\varphi,f,F)$ be a solution of equation \eq{-1} on the interval $I$, $A\in\R$, and $J,K$ be subintervals of $I$ such that $f(J)=f(K)$. If $(\varphi_A,f)$ solves equation \eq{0} on $J$, then it solves equation \eq{0} also on $K$.}

\begin{proof}
Assume that $(\varphi_A,f)$ solves equation \eq{0} on $J$. Then, by \prp{3}, there exists $B\in\R$, such that $F=Af+B$ on $J$. By $f(J)=f(K)$, we get that the same formula holds on $K$, which, in view of \prp{3} again, yields that $(\varphi_A,f)$ solves equation \eq{0} on $K$.
\end{proof}

For a given subinterval $J\subseteq I$, define the sequence of intervals $(J_n)$ by the recursion
\Eq{Jn}{
J_0:=J\qquad\mbox{and}\qquad
J_n:=\Ref(J_{n-1}|J),\qquad(n\in\N).
}
By the definition \eq{Jn}, it is easy to see that $J_n=J$ for all $n\in\N$ provided that $J$ is either empty or is a singleton or equals 
$I$. If $J$ is nonempty then, for the brevity, denote $\inf J_n$ and $\sup J_n$ by $a_n$ and $b_n$ whenever $n\in\N\cup\{0\}$, 
respectively.

\Lem{1}{
Let $J$ be a subinterval of $I$. Then $(J_n)$ is a nondecreasing sequence of intervals contained in $I$ and, provided that $J$ is 
nonempty, for all $n\in\N$, we have
\Eq{ab}{
a_n=\max(\alpha,2a_0-b_{n-1})
\qquad\mbox{and}\qquad
b_n=\min(2b_0-a_{n-1},\beta).
}
Consequently, if $\alpha< a_n$ and $b_n<\beta$ for some $n\in\N\cup\{0\}$ then
\Eq{diam}{
b_k-a_k=(2k+1)(b_0-a_0),\qquad(k\in\{0,1,\dots,n\}).}
Finally, 
if $J\neq I$ is a proper subinterval, then $J=J_0\neq J_1$, furthermore,
if $\alpha=a_n$ or $b_n=\beta$ for some $n\in\N\cup\{0\}$, then $J_{n+k}=J_{n+k+1}$ holds for all $k\in\N$.
}

\begin{proof}
The monotonicity of the sequence $(J_n)$ can be easily shown by induction with respect to $n$. Obviously, $J_0=J\subseteq\Ref(J|J)=J_1$ 
holds. Assume that $J_{n-1}\subseteq J_n$ for some $n\in\N$. Then, by the property (\ref{a}) of \prp{HK}, we have
\Eq{*}{
J_n=\Ref(J_{n-1}|J)\subset \Ref(J_n|J)=J_{n+1}.
}
These imply that $J_{n-1}\subseteq J_n$ for all $n\in\N$. 

Assuming that $J$ is nonempty, the identities in \eq{ab} are easy consequences of \eq{ip} in \prp{HK}. Now suppose that 
$\{a_n, b_n\}\subseteq I$ for some $n\in\N\cup\{0\}$. Then, by the strict monotonicity of $(J_n)$, we have $\{a_k, b_k\}\subseteq I$ for 
all $k\in\{0,1,\dots,n\}$. Using this, we prove \eq{diam} by induction with respect to $k$. If $n=0$ then $k=0$ and the statement is trivial. 
Assume that $n\in\N$. If $k=0$ then the equality in \eq{diam} is trivial. 
Assume that \eq{diam} holds for some non-negative integer $k\leq n-1$, that is, we have 
$b_k-a_k=(2k+1)(b_0-a_0)$. Then, by \eq{ab}, we obtain that
\Eq{*}{
b_{k+1}-a_{k+1}&=2b_0-a_k-(2a_0-b_k)=(b_k-a_k)+2(b_0-a_0)\\
&=(2k+1)(b_0-a_0)+2(b_0-a_0)=\big(2(k+1)+1\big)(b_0-a_0).
}

If $J\neq I$ is a proper subinterval, then $a_0<b_0$ and one of the strict inequalities $\alpha<a_0$ or $b_0<\beta$ holds. In the first case, $a_1<a_0$. Indeed, if this were not true, then by \eq{ab}, $a_1=\max(\alpha,2a_0-b_0)=a_0$.  This imples that $2a_0-b_0=a_0$, hence $a_0=b_0$, a contradiction. In the case $b_0<\beta$, the inequality $b_0<b_1$ follows similarly.

Assume, finally, that there exists $n\in\N$ such that $a_n=\alpha$. Then, by the increasingness on $(J_n)$, for all $k\in\N$, we have
$\alpha\leq a_{n+k}\leq a_n=\alpha$, that is $a_{n+k}=\alpha$. Furthermore, by this and by the second identity in \eq{ab}, it follows that 
$b_{n+k}=\min(2b_0-\alpha,\beta)$ whenever $k\in\N$. Analogously, if $b_n=\beta$ for some $n\in\N$, then $b_{n+k}=\beta$ and $a_{n+k}=\max(\alpha,2a_0-\beta)$ whenever $k\in\N$.
\end{proof}

\Cor{1}{
Using the notations of \lem{1}, we have
\Eq{*}{
J_\infty:=\bigcup\limits_{n\in\N}J_n=
\begin{cases}
J & \mbox{if }J\mbox{ is a singleton},\\[1mm]
\Ref(I|J) & \mbox{otherwise}.
\end{cases}
}}

\begin{proof}
If $J$ is either empty or a singleton or equals $I$, then the statement is trivial. Let assume that $J$ is a proper nonempty subinterval of $I$ 
with $a_0<b_0$. First we show that $J_\infty$ symmetric with respect to $J$, that is, $\Ref(J_\infty|J)=J_\infty$. Indeed,
\Eq{*}{
  \Ref(J_\infty|J)
  =\Big(\bigcup\limits_{n\in\N}J_n\Big|J\Big)
  =\bigcup\limits_{n\in\N}\Ref(J_n|J)
  =\bigcup\limits_{n\in\N}J_{n+1}=J_\infty.
}
This immediatley implies that $J_\infty\subseteq\Ref(I|J)$. To prove the reversed inclusion, let $x\in\Ref(I|J)\setminus J_\infty$ be any point. Then there exist $u\in J$ and $y\in I$ such that $x=2u-y$. We may assume that $x<y$. Obviously, in this case we must have $y\notin J_\infty$. Indeed, if 
$y$ were an element of $J_\infty$, then there would exist $n\in\N$ such that $y\in J_n$. Therefore $x\in\Ref(J_n|J)=J_{n+1}$, which 
contradicts that $x\notin J_\infty$. On the other hand, $u\in J\cap\,]x,y[\,\subseteq J_\infty\cap\,]x,y[$. Consequently, we must have $J_\infty\subseteq\,]x,y[\,$. Using \lem{1}, we obtain that
\Eq{*}{
y-x\geq \diam(J_n)=b_n-a_n=(2n+1)(b_0-a_0)>0,\qquad(n\in\N),
}
which is impossible. Therefore $J_\infty=\Ref(I|J)$.
\end{proof}

\Cor{2}{Let $(\varphi,f,F)$ be a solution of equation \eq{-1} on the interval $I$, $A\in\varphi(I)$ and $J\subseteq\Z_{\varphi_A}$ be a subinterval. Then there exists $B\in\R$ such that $F=Af+B$ holds on $J_\infty$.
}

\begin{proof}
If $\varphi$ is identically $A$ on $J$, then $(\varphi_A,f)$ trivially solves equation \eq{0} on $J$, thus, in view of \prp{3}, $F$ is of the form $Af+B$ on $J=J_0$ for some $B\in\R$. Assume now that this is valid on $J_n$ for some $n\in\N\cup\{0\}$ and let $x\in J_{n+1}$ be any point. Then there exist $u\in J$ and $y\in J_n$ such that $u=\frac{x+y}2$. Then $\varphi(u)=A$ and, by the inductive hypothesis, $F(y)=Af(y)+B$. Thus, applying equation \eq{-1} for $x$ and 
$y$, we get that
\Eq{*}{
A(f(x)-f(y))=F(x)-\big(Af(y)+B\big),
}
which reduces to $F(x)=Af(x)+B$. This means that we have $F=Af+B$ on $J_n$ for arbitrary $n\in\N\cup\{0\}$. Now our statement is an obvious consequence of the monotonicity of the sequence $(J_n)$ and of the definition of the interval $J_\infty$.
\end{proof}

Having these results, we introduce some special subintervals of $I$, which will be useful in the proof of our main theorem.

Suppose that $(\varphi,f,F)$ is a solution of \eq{-1} on $I$. For a given $A\in\varphi(I)$, let $\ubar{I}_A:=\,]\alpha,\beta_A[$ (where $\beta_A\in[\alpha,\beta]$) be the maximal open subinterval $\ubar{I}_A\subseteq I$ such that $(\varphi_A,f)$ solves \eq{0} on $\ubar{I}_A$. Similarly, let $\bar{I}_A:=\,]\alpha_A,\beta[$ (where $\alpha_A\in[\alpha,\beta]$) denote the maximal open subinterval $\bar{I}_A\subseteq I$ such that $(\varphi_A,f)$ solves \eq{0} on $\bar{I}_A$. Finally, define $\ubar{I}$ and $\bar{I}$ by
\Eq{*}{
\ubar{I}:=\bigcup\limits_{A\in\varphi(I)}\ubar{I}_A\qquad\mbox{and}\qquad
\bar{I}:=\bigcup\limits_{A\in\varphi(I)}\bar{I}_A.
}
Obviously, the families $\{\ubar{I}_A\mid A\in\varphi(I)\}$ and $\{\bar{I}_A\mid A\in\varphi(I)\}$ are chains with respect to the inclusion, therefore, $\ubar{I}$ and $\bar{I}$ are (not necessarily nonempty) open subintervals contained in $I$. Then, for some uniquely determined elements $\bar{\alpha},\ubar{\beta}\in[\alpha,\beta]$, we have that $\ubar{I}:=\,]\alpha,\ubar{\beta}[$ and $\bar{I}:=\,]\bar{\alpha},\beta[$.

\Cor{30}{
Let $(\varphi,f,F)$ be a solution of \eq{-1}. Then, for all $A\in\varphi(I)$ and for all subinterval $J\subseteq\Z_{\varphi_A}$, we have $J_\infty\subseteq\ubar{I}$ or $J_\infty\subseteq\bar{I}$.
}

\begin{proof}
Let $A\in\varphi(I)$ and $J\subseteq\Z_{\varphi_A}$ be any subinterval. Then, by \cor{2}, there exists $B\in\R$, such that $F=Af+B$ on $J_\infty$. This, in view of \prp{3}, is equivalent to that the pair $(\varphi_A,f)$ solves the equation \eq{0} on the interval $J_\infty$. Hence at least one of the inclusions $J_\infty\subseteq\ubar{I}$ or $J_\infty\subseteq\bar{I}$ must hold.
\end{proof}

The following lemmas are about some useful connections between the structure of the sets $\ubar I$ and $\bar I$, and the behavior of the members of the solutions triplets of equation \eq{-1}.

\Lem{4}{
Let $(\varphi,f,F)$ be a solution of \eq{-1}, $A$ and $A'$ be distinct points in $\varphi(I)$, furthermore $J$ and $K$ be intervals such that $(\varphi_A,f)$ and $(\varphi_{A'},f)$ solves \eq{0} on $J$ and $K$, respectively. Then $f$ is constant on $J\cap K$.
}

\begin{proof}
If $J\cap K$ is empty, then the statement is trivial. If not, then we have
\Eq{*}{
\Big(\varphi\Big(\frac{x+y}{2}\Big)-A\Big)(f(x)-f(y))=0\qquad\mbox{and}\qquad
\Big(\varphi\Big(\frac{x+y}{2}\Big)-A'\Big)(f(x)-f(y))=0
}
for all $x,y\in J\cap K$. Subtracting these equations side by side, we get $(A-A')(f(x)-f(y))=0$, which, by our assumption $A\neq A'$, implies that $f(x)=f(y)$.
\end{proof}

\Lem{5}{
Let $(\varphi,f,F)$ be a solution of \eq{-1}. Then $f$ is constant on $\ubar{I}$ (resp. on $\bar{I}$) or there uniquely exists $A\in\varphi(I)$ such that $\ubar{I}=\ubar{I}_A$ (resp. $\bar{I}=\bar{I}_A$).
}

\begin{proof}
We prove the statement only for $\ubar{I}$. To avoid the trivial case, we may assume that $\ubar{I}$ is nonempty. Assume that there is no $A\in\varphi(I)$ such that $\ubar{I}=\ubar{I}_A$. Then it is enough to show that $f$ is constant on any member of the chain-union in the definition of $\ubar{I}$. To do this, let $A\in\varphi(I)$ be arbitrary and $x\in\ubar{I}\setminus\ubar{I}_A$. Then, there exists $A'\in\varphi(I)$ different from $A$ such that $x\in\ubar{I}_{A'}$. Hence $\ubar{I}_A\subseteq\ubar{I}_{A'}$ and, in view of \lem{4}, the function $f$ is constant on $\ubar{I}_A\cap\ubar{I}_{A'}=\ubar{I}_A$.
\end{proof}

\Lem{cp}{Let $(\varphi,f,F)$ be a solution of \eq{-1}. Then, for all constants $\lambda,\mu\in\R$ with $\lambda\neq\mu$, the function $\varphi$ is constant on the set $\frac12(\Z_{f-\lambda}+\Z_{f-\mu})$.
}

\begin{proof}
If at least one of the sets $\Z_{f-\lambda}$ and $\Z_{f-\mu}$ is empty, then $\frac12(\Z_{f-\lambda}+\Z_{f-\mu})$ is also empty, thus, in this case, the statement is trivially true. Assume therefore that $\Z_{f-\lambda}$ and $\Z_{f-\mu}$ are nonempty. Equation \eq{-1} implies that $F$ also must be constant both on the sets $\Z_{f-\lambda}$ and $\Z_{f-\mu}$. Thus there exists not necessarily different numbers $C_\lambda,C_\mu\in\R$ such that $F(u)=C_\lambda$ and $F(v)=C_\mu$ whenever $u\in\Z_{f-\lambda}$ and $v\in\Z_{f-\mu}$. Now, applying equation \eq{-1} for any points $x\in\Z_{f-\lambda}$ and $y\in\Z_{f-\mu}$, we obtain that
\Eq{*}{
\varphi\Big(\frac{x+y}2\Big)
=\frac{F(x)-F(y)}{f(x)-f(y)}
=\frac{C_\lambda-C_\mu}{\lambda-\mu}.
}
Thus $\varphi$ is indeed constant on $\frac12(\Z_{f-\lambda}+\Z_{f-\mu})$.
\end{proof}

\Lem{31}{
Let $(\varphi,f,F)$ be a solution of \eq{-1}, furthermore assume that the closed subinterval $S:=I\setminus(\ubar{I}\cup\bar{I})$ has a nonempty interior and that it is different from $I$. Then the following statements hold.
\begin{enumerate}[(i)]\itemsep=1mm
\item\label{i} The function $\varphi$ is locally non-constant on $S$.
\item\label{ii} If there exists $\gamma\in[\alpha,\ubar\beta[\,$ such that $f$ is constant on the interval $]\gamma,\ubar\beta[\,$, then $f$ must be constant on $]\gamma,\bar\alpha[$. Similarly,
if there is $\gamma\in\,]\bar\alpha,\beta]$ such that $f$ is constant on $]\bar\alpha,\gamma[\,$, then $f$ must be constant on $]\ubar\beta,\gamma[$.
\end{enumerate}
}

\begin{proof}
To prove \emph{(i)}, indirectly, assume that there exists $A\in\varphi(I)$ such that the interior of the intersection $S\cap\Z_{\varphi_A}$ is nonempty. If 
$J\subseteq S\cap\Z_{\varphi_A}$ is any nonempty open subinterval, then, based on \cor{1}, $J\subseteq J_\infty$ holds. Furthermore, by \cor{30}, we obtain that
$J\subseteq S\cap J_{\infty}\subseteq S\cap(\ubar{I}\cup\bar{I})$, where the last intersection is empty, contradicting that $J$ was a proper open subinterval.
	
Now we prove the first assertion of statement \emph{(ii)}. Assume that there exists $\alpha\leq\gamma<\ubar\beta$ and $\lambda\in\R$ such that $f$ is identically $\lambda$ on the interval $]\gamma,\ubar\beta[\,$. Let $y$ be any point of the interior of $S$. Then there exists $x\in\,]\gamma,\ubar\beta[$ such that $\frac{x+y}{2}$ belongs also to the interior of $S$. If $f(y)\neq\lambda$, then, based on \lem{cp}, $\varphi$ must be constant on $\frac12(\{y\}+[x,\ubar\beta[)\subseteq S$, which, in view of the statement \emph{(i)}, is impossible. This shows that $f(y)$ must equal $\lambda$.

A similar argument yields that $f(\ubar\beta)=\lambda$ also holds. Otherwise, again by \lem{cp}, $\varphi$ would be constant on $\frac12(\ubar\beta+\,]\ubar\beta,\bar\alpha[)\subseteq S$, which contradicts \emph{(i)} again. Thus, $f$ is indeed constant on $]\gamma,\bar\alpha[\,$.

The proof of the second assertion of \emph{(ii)} can be treated similarly.
\end{proof}

\Thm{2}{
Let $\varphi:I\to\R$ be a function with closed level sets in $I$ and assume that there exists a level set of $\varphi$ whose interior is nonempty. Then the triplet of functions $(\varphi,f,F)$ solves functional equation \eq{-1} if and only if 
there exist constants $A\in\varphi(I)$ and $B\in\R$ such that $(\varphi_A,f)$ solves \eq{0} on $I$ and $F=Af+B$ on $I$.}

\begin{proof}
The sufficiency is trivial. Indeed, assuming there exist $A\in\varphi(I)$ and $B\in\R$ such that $(\varphi_A,f)$ solves \eq{0} on $I$ and $F=Af+B$ on $I$ and calculating both sides of equation the \eq{-1}, one can conclude that $(\varphi,f,F)$ solves functional equation \eq{-1}. We note that the argument in this direction does not need any condition concerning the level sets of $\varphi$.

Now we turn to the necessity, that is, assume that $(\varphi,f,F)$ solves functional equation \eq{-1}. In view of \prp{3}, it is enough to show that there exists $A\in\varphi(I)$, such that $(\varphi_A,f)$ solves equation \eq{0} on the interval $I$. This is trivially fulfilled provided that at least one of the functions $f$ or $\varphi$ is constant on $I$. Therefore, in the sequel, we suppose that $f$ and $\varphi$ are non-constant functions.

By our assumption, $\varphi$ has at least one level set with a nonempty interior, which, based on \cor{30}, means that at least one of the sets $\ubar{I}$ and $\bar{I}$ is nonempty. Without loss of the generality, we may assume that $\ubar{I}\neq\emptyset$. If $\ubar{I}=I$, then we are done. Hence, in the rest of the proof we may also assume that $\ubar{I}$ is a proper subinterval of $I$.

For the brevity, denote the interval $I\setminus(\ubar{I}\cup\bar{I})$ by $S$. We are going to prove that $S$ must be empty. More precisely, first we show that its interior is empty, then, secondly, that it cannot be a singleton.

\emph{Part 1.} Indirectly assume that the interior of $S$ is nonempty. Observe, that, by the statement \emph{(ii)} of \lem{31}, the function $f$ cannot be constant on $\ubar I$. Indeed, if $f$ were constant on $\ubar I=\,]\alpha,\ubar\beta[\,$, then it would also be constant on $]\alpha,\bar{\alpha}[\,=\ubar I\cup(S\setminus\{\bar\alpha\})$. Then, for any $A$, the pair $(\varphi_A,f)$ satisfies \eq{0} on $]\alpha,\bar{\alpha}[$, which means that $\ubar I$ would be strictly expandable, contradicting its definition. Therefore, based on \lem{5}, there uniquely exists $A\in\varphi(I)$ such that $\ubar I=\ubar I_A$, that is, the pair $(\varphi_A,f)$ solves equation \eq{0} on the subinterval $\ubar I$. Applying \thm{1} for the interval $\ubar I$, there exists a closed subinterval $K\subseteq\ubar I$ such that $J:=\frac12(\ubar I+K)\subseteq\Z_{\varphi_A}$ holds. 

Using the notations of \thm{1}, we show that $K^*\cap \ubar I$ must be empty. If not, then, in view of the statement \emph{(ii)} of \lem{31}, we have that $f$ is constant on $K^*\cap\,]\alpha,\bar\alpha[\,$. This implies that $f(\ubar I)=f(]\alpha,\bar\alpha[)$.  
Using that the pair $(\varphi_A,f)$ solves equation \eq{0} on $\ubar I$, by \cor{ext}, it follows that $(\varphi_A,f)$ also solves equation \eq{0} on $]\alpha,\bar\alpha[$. This contradicts the maximality property of $\ubar I$ again.

Therefore, $K^*\cap \ubar I$ is empty, or equivalently, we have that $\sup J:=\ubar\beta$. By \lem{1}, the interval $J_1$ and hence $J_\infty$ contains the point $\ubar\beta$ in its interior. On the other hand, \cor{30} implies that $J_\infty\subseteq \ubar I$ or $J_\infty\subseteq \bar I$. Both inclusions are obviously impossible. This final contradiction means that the interior of $S$ must be indeed empty or, equivalently, the closed interval $S$ is at most a singleton.

\emph{Part 2.} In this part we are going to show that $S$ is neither a singleton. The proof of this part also goes indirectly, that is, let us assume that $S=\{p\}\subseteq I$. Then both of the open intervals $\ubar I$ and $\bar I$ are nonempty. Furthermore, there is no one-sided open neighborhood of $p$ contained in $I$, where $\varphi$ were constant. Otherwise, if there exists $r>0$ such that $\varphi$ is constant on the left neighborhood $J:=\,]p-r,p]\subseteq I$ of $p$. Then the interval $J_\infty$ contains $p$ in its interior, which, because of \cor{30}, is impossible. A similar argument shows that $\varphi$ cannot be constant on any right neighborhood of $p$.

Therefore, in view of \lem{5} and \thm{1}, there exist nonempty open subintervals $L, R\subseteq I$ with $\sup L=\inf R=p$ and real constants $\lambda_L,\lambda_R\in\R$ such that $f|_{L}\equiv\lambda_L$ and $f|_{R}\equiv\lambda_R$. We claim that $\lambda_L$ and $\lambda_R$ are necessarily equal to $f(p)$. If not, then, by \lem{cp}, we get that $\varphi$ is constant on the interval $\frac12(L+{p})$ and $\frac12({p}+R)$, respectively. These sets are left and right neighborhoods of $p$, which, in view of the previous argument, yields a contradiction.

Thus the restriction of $f$ to the interval $L\cup\{p\}\cup R$ is a constant function. Then we have $f(\ubar I)=f(\ubar I\cup R)$. Since $(\varphi_A,f)$ solves equation \eq{0} on $\ubar I$, by \cor{ext}, it follows that $(\varphi_A,f)$ also solves equation \eq{0} on $]\ubar I\cup R[$. This again contradicts the maximality property of $\ubar I$.

Consequently, $S$ is empty or, equivalently, the open interval $\ubar I\cap\bar I$ is nonempty. Particularly, we get that $\bar I$ is also nonempty. To avoid the trivial case, it can be supposed that $\bar I$ is different from $I$. The function $f$ is non-constant on $I$, therefore we may assume that $f$ is non-constant on $\ubar I$. Then, by \lem{5}, there uniquely exists $A\in\varphi(I)$ such that $\ubar I=\ubar I_A$. If $f$ is constant on $\bar I$, then $f(\ubar I)=f(\ubar I\cup\bar I)=f(I)$. Thus, in view of \cor{ext}, the pair $(\varphi_A,f)$ solves equation \eq{0} on the entire domain $I$, hence in this case we are done. Therefore in the sequel we may assume also that $f$ is non-constant on $\bar I$. By \lem{5}, there uniquely exists $A'\in\varphi(I)$ such that $\bar I=\bar I_{A'}$. Now, obviously, we have two possibilities, namely either $A\neq A'$ or $A=A'$.

Assume that $A\neq A'$ holds. Then, by \lem{4}, there exists a real constant $\lambda$ such that $f$ is identically $\lambda$ on the intersection $M:=\ubar I\cap\bar I$. Then, because of the maximality property of $\ubar I$ and $\bar I$, there is no $r>0$ such that $f$ takes $\lambda$ also on the open subinterval $(M\,+\,]-r,r[\,)\cap I$. Consequently, for all $n\in\N$, there exists $y_n\in\,]\ubar\beta,\ubar\beta+\frac1n[\,\cap\, I$ such that $f(y_n)\neq\lambda$. Then, by \lem{cp}, the function $\varphi$ must be constant on the interval $\frac12\big((\ubar I\cap\bar I)+y_n\big)$ for all $n\in\N$. The intersection of the family $\big\{\frac12\big((\ubar I\cap\bar I)+y_n\big)\mid n\in\N\big\}$ is nonempty, which implies that $\varphi$ is constant on the union $\bigcup_{n\in\N}\frac12\big((\ubar I\cap\bar I)+y_n\big)\supseteq\,\big]\frac{\bar\alpha+\ubar\beta}{2},\ubar\beta\big[\,$.
A similar argument shows that $\varphi$ must be also constant on the interval $\big]\bar\alpha,\frac{\bar\alpha+\ubar\beta}{2}\big[\,$. The level sets of $\varphi$ are closed, therefore $\varphi$ takes the same value on $J:=]\bar\alpha,\ubar\beta[\,$. Then $J_\infty$ contains the points $\bar\alpha$ and $\ubar\beta$ in its interior, which contradicts the statement of \cor{30}.

Finally, we obtained that $A=A'$ must hold. Then the pair $(\varphi_A,f)$ trivially solves \eq{0} on the sets $\ubar I$ and $\bar I$. By \prp{3}, this means that there exist $B,B'\in\R$ such that $F=Af+B$ on $\ubar I$ and $F=Af+B'$ on $\bar I$. The nonemptyness of the intersection $\ubar I_A\cap\bar I_{A}$ yields that $B=B'$. Consequently, $F=Af+B$ on $I$. This, based again on \prp{3}, equivalent to that $(\varphi_A,f)$ solves \eq{0} on the entire interval $I$, which finishes the proof.
\end{proof}

\Cor{2+}{Let $\varphi:I\to\R$ be a function with closed level sets in $I$ and assume that there exists a level set of $\varphi$ whose interior is nonempty. Then the triplet of functions $(\varphi,f,F)$ solves functional equation \eq{-1} if and only if either $f$ is constant on $I$ and $\varphi,F:I\to\R$ are arbitrary functions or there exist $A\in\varphi(I)$, $B\in\R$, a nonempty closed interval $K\subseteq I$, and constants $\lambda_*,\lambda^*\in\R$ such that the inclusions
\Eq{Sol2}{
K_*\subseteq\Z_{f-\lambda_*},
\qquad
K^*\subseteq\Z_{f-\lambda^*}
\qquad\mbox{and}\qquad
\tfrac12(I+K)\subseteq\Z_{\varphi-A}
}
are satisfied, furthermore we have $F=Af+B$ on $I$.
}
\section{Solutions of \eq{-1} assuming that $\varphi$ is continuous}

The following statement summarizes what we have established in the previous sections and establishes the main result of our paper.

\Thm{summ}{
Let $\varphi:I\to\R$ be a continuous function. Then the triplet of functions $(\varphi,f,F):I\to\R^3$ solves equation \eq{-1} on $I$ if and only if either
\begin{enumerate}[(i)]\itemsep=1mm
\item $f$ and $F$ are constant functions and $\varphi$ is an arbitrary function on $I$, or
\item there exist $A\in\varphi(I)$, a nonempty closed interval $K\subseteq I$ and constants $\lambda_*,\lambda^*,\mu\in\R$ such that the inclusions of 
\eq{Sol2} hold and $F=Af+\mu$ on $I$, or
\item there exist constants $A,B,C,D\in\R$ with $AD\neq BC$ and $\gamma,\mu,\lambda\in\R$ such that, for all $x\in I$, one of the following possibilities hold.
\Eq{neg}{
\begin{array}{c}
 F(x)=-C\cos(\sqrt{-\gamma}x)+D\sin(\sqrt{-\gamma}x)+\mu, \\[1mm]
 f(x)=-A\cos(\sqrt{-\gamma}x)+B\sin(\sqrt{-\gamma}x)+\lambda,
\end{array}
\quad 
\varphi(x)=\frac{C\sin(\sqrt{-\gamma}x)+D\cos(\sqrt{-\gamma}x)}{A\sin(\sqrt{-\gamma}x)+B\cos(\sqrt{-\gamma}x)}
\quad\mbox{if }\gamma<0,
}
\Eq{null}{\,\,
\begin{array}{c}
 F(x)=\frac12Cx^2+Dx+\mu,\\[1mm]
 f(x)=\frac12Ax^2+Bx+\lambda,
\end{array}
\quad
\varphi(x)=\frac{Cx+D}{Ax+B}
\qquad\mbox{if }\gamma=0,
}
\Eq{poz}{
\begin{array}{c}
 F(x)=C\cosh(\sqrt{\gamma}x)+D\sinh(\sqrt{\gamma}x)+\mu,\\[1mm]
 f(x)=A\cosh(\sqrt{\gamma}x)+B\sinh(\sqrt{\gamma}x)+\lambda,
\end{array}
\quad
\varphi(x)=
\frac{C\sinh(\sqrt{\gamma}x)+D\cosh(\sqrt{\gamma}x)}{A\sinh(\sqrt{\gamma}x)+B\cosh(\sqrt{\gamma}x)}
\qquad\mbox{if }\gamma>0.
}
\end{enumerate}
}

It is worth noticing that the condition (iii) of the theorem could be formulated in the following equivalent way:
\begin{enumerate}
 \item[(iii)'] $f$ and $F$ are differentiable functions, $f'$ is nowhere zero on $I$ and there exists a constant $\gamma\in\R$ such that $f'$ and $F'$ are linearly independent solutions of the second-order linear homogeneous differential equation $Y''=\gamma Y$ and $\varphi=F'/f'$ on $I$.
\end{enumerate}

\begin{proof} The continuity of $\varphi$ implies that, for all $E\in\R$, the level set $\Z_{\varphi-E}$ is closed in $I$.

It is obvious that if $f$ and $F$ are constants and $\varphi$ is an arbitrary function on $I$, then $(\varphi,f,F):I\to\R^3$ solves equation \eq{-1} on $I$. If $F$ is constant on $I$, then $(\varphi,f)$ solves \eq{0}. Therefore, by \thm{1}, either $f$ is constant on $I$ or there exist a nonempty closed interval $K\subseteq I$ and constants $\lambda_*,\lambda^*\in\R$ such that the inclusions in \eq{Sol0} hold. Thus the alternative \emph{(ii)} of our theorem is valid with $A=0$. Therefore, in the sequel, we may assume that $F$ is non-constant on $I$. Then $f$ is also non-constant.

If there exists a level set of $\varphi$ whose interior is nonempty, then, in view of \cor{2+}, there exist $A\in\varphi(I)$, a nonempty closed interval $K\subseteq I$ and constants $\lambda_*,\lambda^*,\mu\in\R$ such that the inclusions \eq{Sol2} hold, furthermore $F=Af+\mu$ on $I$. This shows, that the alternative \emph{(ii)} is valid in this case.

Finally, we consider the case when $\varphi$ is locally non-constant and the functions $f$ and $F$ are non-constant on $I$. We are going to show that in this setting the alternative \emph{(iii)} is valid. First we deal with the necessity of \emph{(iii)}. Assume that the triplet of functions $(\varphi,f,F):I\to\R^3$ solves equation \eq{-1} on $I$. Then, in view of \thm{0}, the function $f$ is infinitely many times differentiable and $(\varphi,f')$ satisfies the functional equation \eq{f'} on $I$. Consequently, one of the possibilities \emph{(i)} or \emph{(ii)} of \thm{m1} holds concerning the pair $(\varphi,f')$.

If the alternative \emph{(i)} of \thm{m1} were valid, then there would exist a subinterval $J\subseteq I$ such that $f'$ is zero on $I\setminus J$ and $\varphi$ is constant on $\frac12(J+I)$. Using that $\varphi$ is locally non-constant on $I$, we could conclude that $J$ must be empty, which would yield that $f$, and hence $F$, is constant on $I$ contradicting our assumptions. Therefore, we must have the alternative \emph{(ii)} of \thm{m1}. Then there exist real constants $a,b,c,d,\gamma\in\R$ with the property $ad\neq bc$ such that, depending on the sign of $\gamma$, the functions $\varphi$ and $f'$ are of the forms \eq{m1a}, \eq{m1b} or \eq{m1c} listed in \thm{m1}. To compute $f$ and $F$, we distinguish three main cases.

\emph{Case $\gamma<0$.} Then, in view of (3) of \thm{m1}, we have
\Eq{*}{
\varphi(x)=
\frac{c\sin(\sqrt{-\gamma}x)+d\cos(\sqrt{-\gamma}x)}{a\sin(\sqrt{-\gamma}x)+b\cos(\sqrt{-\gamma}x)}\qquad\mbox{and}\qquad
f'(x)=a\sin(\sqrt{-\gamma}x)+b\cos(\sqrt{-\gamma}x)
}
for all $x\in I$. Therefore, there exists $\lambda\in\R$ such that
\Eq{*}{
f(x)=\frac{1}{\sqrt{-\gamma}}(-a\cos(\sqrt{-\gamma}x)+b\sin(\sqrt{-\gamma}x))+\lambda,\qquad(x\in I).
}
Replacing $a,b,c$, and $d$ by $A\sqrt{-\gamma}$, $B\sqrt{-\gamma}$, $C\sqrt{-\gamma}$, and $D\sqrt{-\gamma}$, respectively, we can see that $f$ and $\varphi$ are of the form stated in \eq{neg}. To conclude the proof, we compute $\varphi\big(\frac{x+y}2\big)(f(x)-f(y))$ for $x,y\in I$. To simplify the calculation, denote the points $\sqrt{-\gamma}x$ and $\sqrt{-\gamma}y$ by $u$ and $v$, respectively. Then
\Eq{*}{
  f(x)-f(y)
  &=(-A\cos(u)+B\sin(u))-(-A\cos(v)+B\sin(v))
  =2\sin\tfrac{u-v}{2}\big(A\sin\tfrac{u+v}{2}+B\cos\tfrac{u+v}{2}\big).
}
Therefore,
\Eq{id-}{
 \varphi\Big(\frac{x+y}2\Big)(f(x)-f(y))
 &=\frac{C\sin\frac{u+v}{2}+D\cos\frac{u+v}{2}}{A\sin\frac{u+v}{2}+B\cos\frac{u+v}{2}}\cdot
 2\sin\tfrac{u-v}{2}\big(A\sin\tfrac{u+v}{2}+B\cos\tfrac{u+v}{2}\big)\\
 &=2\sin\tfrac{u-v}{2}\big(C\sin\tfrac{u+v}{2}+D\cos\tfrac{u+v}{2}\big)\\
 &=(-C\cos(u)+D\sin(u))-(-C\cos(v)+D\sin(v))\\
 &=(-C\cos(\sqrt{-\gamma}x)+D\sin(\sqrt{-\gamma}x))-(-C\cos(\sqrt{-\gamma}y)+D\sin(\sqrt{-\gamma}y)).
}
Using that $(\varphi,f,F)$ solves \eq{-1}, the above identity yields, for all $x,y\in I$, that
\Eq{*}{
  F(x)-F(y)=(-C\cos(\sqrt{-\gamma}x)+D\sin(\sqrt{-\gamma}x))-(-C\cos(\sqrt{-\gamma}y)+D\sin(\sqrt{-\gamma}y)),
}
hence the mapping $x\mapsto F(x)+C\cos(\sqrt{-\gamma}x)-D\sin(\sqrt{-\gamma}x)$ is constant on $I$. This shows that $F$ is of the form \eq{neg} for some real number $\mu$.

\emph{Case $\gamma=0$.} Then, by (2) of \thm{m1}, for all $x\in I$, we have
\Eq{*}{
\varphi(x)=\frac{cx+d}{ax+b}\qquad\mbox{and}\qquad
f'(x)=ax+b.
}
As before, we immediately get that there exists $\lambda\in\R$ such that
\Eq{*}{
f(x)=\frac{a}2x^2+bx+\lambda,\qquad(x\in I).
}
Now, replacing $a,b,c,$ and $d$ by $A,B,C,$ and $D$, respectively, we get, for all $x,y\in I$, that
\Eq{id0}{
\varphi\Big(\frac{x+y}2\Big)(f(x)-f(y))
&=\frac{C\frac{x+y}{2}+D}{A\frac{x+y}{2}+B}\Big(\frac12Ax^2+Bx-\frac12Ay^2-By\Big)
\\&=\frac{C\frac{x+y}{2}+D}{A\frac{x+y}{2}+B}\Big(A\frac{x+y}2+B\Big)(x-y)
=\Big(C\frac{x+y}2+D\Big)(x-y)\\
&=\frac12Cx^2+Dx-\Big(\frac12Cy^2+Dy\Big).
}
Using that $(\varphi,f,F)$ solves \eq{-1}, the above identity yields, for all $x,y\in I$, that we must have
\Eq{*}{
F(x)-F(y)=\frac12Cx^2+Dx-\Big(\frac12Cy^2+Dy\Big),
}
hence the mapping $x\mapsto F(x)-\frac12Cx^2-Dx$ is constant on $I$. This shows that $F$ is of the form \eq{null} for some real number $\mu$.

\emph{Case $\gamma>0$.} Based on (3) of \thm{m1}, in this case we have
\Eq{*}{
\varphi(x)=
\frac{c\sinh(\sqrt{\gamma}x)+d\cosh(\sqrt{\gamma}x)}{a\sinh(\sqrt{\gamma}x)+b\cosh(\sqrt{\gamma}x)}\qquad\mbox{and}\qquad
f'(x)=a\sinh(\sqrt{\gamma}x)+b\cosh(\sqrt{\gamma}x)
}
for all $x\in I$, thus there exists $\lambda\in\R$ such that
\Eq{*}{
f(x)=\frac{1}{\sqrt{\gamma}}(a\cosh(\sqrt{\gamma}x)+b\sinh(\sqrt{\gamma}x))+\lambda,\qquad(x\in I).
}
Substituting $A\sqrt{\gamma}$, $B\sqrt{\gamma}$, $C\sqrt{\gamma}$, and $D\sqrt{\gamma}$ instead of the constants $a,b,c$, and $d$, respectively, we can see that $f$ and $\varphi$ are of the form stated in \eq{poz}. Now we compute $\varphi\big(\frac{x+y}2\big)(f(x)-f(y))$ for $x,y\in I$. For the brevity, denote the elements $\sqrt{\gamma}x$ and $\sqrt{\gamma}y$ by $u$ and $v$, respectively. Then,
\Eq{*}{
  f(x)-f(y)
  &=(A\cosh(u)+B\sinh(u))-(A\cosh(v)+B\sinh(v))
  =2\sinh\tfrac{u-v}{2}\big(A\sinh\tfrac{u+v}{2}+B\cosh\tfrac{u+v}{2}\big).
}
Therefore, for all $x,y\in I$, we have
\Eq{id+}{
 \varphi\Big(\frac{x+y}2\Big)(f(x)-f(y))
 &=\frac{C\sinh\frac{u+v}{2}+D\cosh\frac{u+v}{2}}{A\sinh\frac{u+v}{2}+B\cosh\frac{u+v}{2}}\cdot
 2\sinh\tfrac{u-v}{2}\big(A\sinh\tfrac{u+v}{2}+B\cosh\tfrac{u+v}{2}\big)\\
 &=2\sinh\tfrac{u-v}{2}\big(C\sinh\tfrac{u+v}{2}+D\cosh\tfrac{u+v}{2}\big)\\
 &=(C\cosh(u)+D\sinh(u))-(C\cosh(v)+D\sinh(v))\\
 &=(C\cosh(\sqrt{\gamma}x)+D\sinh(\sqrt{\gamma}x))-(C\cosh(\sqrt{\gamma}y)+D\sinh(\sqrt{\gamma}y)).
}
Using that $(\varphi,f,F)$ solves \eq{-1}, the above identity yields, for all $x,y\in I$, that
\Eq{*}{
  F(x)-F(y)=(C\cosh(\sqrt{\gamma}x)+D\sinh(\sqrt{\gamma}x))-(C\cosh(\sqrt{\gamma}y)+D\sinh(\sqrt{\gamma}y)),
}
hence the function $x\mapsto F(x)-C\cosh(\sqrt{\gamma}x)-D\sinh(\sqrt{\gamma}x)$ is constant on $I$. This verifies that $F$ is of the form \eq{poz} for some $\mu\in\R$.

To show the sufficiency of alternative \emph{(iii)}, observe that in each of the above cases, the identities \eq{id-}, \eq{id0}, and \eq{id+} are satisfied, respectively. Therefore, the triplet $(\varphi,f,F)$ satisfies functional equation \eq{-1}.
\end{proof}

\Rem{A}{One can easily see that the triplets $\varphi,f,F)$ described in alternatives (i) and (ii) are solutions of \eq{-1} without assuming the continuity of $\varphi$. It remains an open problem if the continuity assumption about $\varphi$ can be omitted from the formulation of \thm{summ}.}

\section{An application}
In this section, using our main result, we solve the equality problem of two-variable Cauchy means and two-variable quasi-arithmetic means solely under the conditions needed for their definitions.

\Thm{CA}{
Let $J\subseteq\R$ be an open subinterval. Assume that $G,H:J\to\R$ are differentiable functions such that $0\notin H'(J)$ and $\frac{G'}{H'}$ is invertible on $J$, and that $\Phi:J\to\R$ is continuous and strictly monotone. Then the functional equation
\Eq{eqCA}{
\C_{G,H}(x,y)=\A_\Phi(x,y),\qquad(x,y\in J)}
holds if and only if $\Phi$ is differentiable with a nonvanishing first derivative and there exist constants $A,B,C,D\in\R$ with $AD\neq BC$ and $\mu,\lambda,\gamma\in\R$ such that
\Eq{aff}{
\begin{pmatrix}
G \\ H
\end{pmatrix}=
\begin{pmatrix}
A & B \\ C & D
\end{pmatrix}
\begin{pmatrix}
\psi_1 \\ \psi_2
\end{pmatrix}
+
\begin{pmatrix}
\mu \\ \lambda
\end{pmatrix}
}
holds on $J$, where
\Eq{cases}{
(\psi_1(x), \psi_2(x))=
\begin{cases}
\big(\cos(\sqrt{-\gamma}\Phi(x)), \sin(\sqrt{-\gamma}\Phi(x))\big)
& \mbox{if }\gamma<0, \\[1mm]
\big(\Phi^2(x), \Phi(x)\big) & \mbox{if }\gamma=0, \\[1mm]
\big(\cosh(\sqrt{\gamma}\Phi(x)), \sinh(\sqrt{\gamma}\Phi(x))\big)
 & \mbox{if }\gamma>0,
\end{cases}\qquad\qquad(x\in J).
}}
\begin{proof}
A simple calculation yields that, the triplet $(G,H,\Phi)$ solves functional equation \eq{eqCA} on $J$ if and only if
the triplet of functions $(\varphi,f,F):=\big(\frac{G'}{H'}\circ\Phi^{-1},H\circ\Phi^{-1},G\circ\Phi^{-1}\big)$ solves the functional equation \eq{-1} on $I:=\Phi(J)$.
In view of this connection, the sufficiency part of the statement is obvious. Thus we can turn to the proof of the necessity part, that is, assume that \eq{eqCA} is satisfied. 

Due to the properties of $\Phi$, the set $\Phi(J)$ is an open subinterval of $\R$. By our assumption, the function $\frac{G'}{H'}$ is invertible on $J$, which, in view of Remark 1 of the paper \cite{Mat04c}, implies that it is continuous and strictly monotone. This and the assumptions concerning $\Phi$ imply that $\varphi=\frac{G'}{H'}\circ\Phi^{-1}$ must be also continuous, hence the conditions of \thm{summ} are satisfied. Consequently, for the members of the triplet $(\varphi,f,F)$, we have one the possibilities \emph{(i)}, \emph{(ii)} or \emph{(iii)} listed in \thm{summ}.

Assume that there exists a nonempty open subinterval of $I$, where $f$ is constant. Then, by the continuity and strict monotonicity of $\Phi$, the preimage of this interval is an open interval again, contained in $J$, where $H$ is turned to be also constant. This contradicts the assumption $0\notin H'(J)$. This means that $f$ must be locally non-constant on $I$, which means that the case \emph{(i)} is excluded.

A similar argument shows that the case \emph{(ii)} of \thm{summ} is also impossible. Indeed, assume indirectly that we have \emph{(ii)}. According to the previous part, the set $K$ must coincide with $I$, that is, $\varphi$ must be constant on $\frac12(I+K)=I$. Due to the properties of $\Phi$, this yields that there exists a nonempty open subinterval of $J$, where the ratio $\frac{G'}{H'}$ is constant, which contradicts that it is invertible.

Consequently, there exist real constants $A,B,C$, and $D$ with the property $AD\neq BC$ and $\gamma\in\R$ such that one of the alternatives \eq{neg}, \eq{null} or \eq{poz} listed in the case \emph{(iii)} of \thm{summ} must hold.

Assume that $\gamma<0$ holds. Then, using the definitions of the functions $\varphi,f$ and $F$, for all $u\in I=\Phi(J)$, we have
\Eq{*}{
F(u)=G\big(\Phi^{-1}(u)\big)=-A\cos(\sqrt{-\gamma}u)+B\sin(\sqrt{-\gamma}u)+\mu}
and
\Eq{*}{
f(u)=H\big(\Phi^{-1}(u)\big)=-C\cos(\sqrt{-\gamma}u)+D\sin(\sqrt{-\gamma}u)+\lambda.}
Writing $x$ instead of $\Phi^{-1}(u)$, defining $\psi_1(x):=\cos(\sqrt{-\gamma}\Phi(x))$ and $\psi_2(x):=\sin(\sqrt{-\gamma}\Phi(x))$ for $x\in J$, and, finally, writing $A$ and $C$ instead of $-A$ and $-C$, we get that the pair $(G,H)$ is of the form written in \eq{aff}. The proof in the cases $\gamma=0$ and $\gamma>0$ goes analogously.

Now we are able to prove that $\Phi$ is differentiable with a nonvanishing first derivative. By the Rolle Mean Value Theorem, $H$ has to be invertible on $I$, hence we have $\Phi^{-1}=H^{-1}\circ f$. In view of \thm{summ}, the function $f$ is (continuously) differentiable and, a short calculation yields that its derivative is either positive or negative on its domain. Hence, taking the inverse of both sides of this last identity, we get that $\Phi=f^{-1}\circ H$, which shows that $\Phi$ is differentiable. Finally, the assumption $0\notin H'(J)$ provides that $\Phi'$ does not vanish on $J$.
\end{proof}



\end{document}